\documentstyle[12pt]{article}
\begin{document}

\voffset=-1truecm

\setlength{\textheight}{21cm}

\def\dsp{\displaystyle}
\def\Rr{{bf R}}
\def\Zz{bf Z}
\def\Nn{bf N}
\def\get{\hbox{{\goth g}$^*$}}
\def\g{\gamma}
\def\om{\omega}
\def\r{\rho}
\def\a{\alpha}
\def\s{\sigma}
\def\vfi{\varphi}
\def\l{\lambda}
\def\implique{\Rightarrow}
\def\o{{\circ}}
\def\Diff{\hbox{\rm Diff}}
\def\S1{\hbox{\rm S$^1$}}
\def\Hom{\hbox{\rm Hom}}
\def\Vect{\hbox{\rm Vect}}
\def\const{\hbox{\rm const}}
\def\ad{\hbox{\hbox{\rm ad}}}
\def\semid{\hbox{\bb o}}
\def\blanc{\hbox{\ \ }}

\def\pds#1,#2{\langle #1\mid #2\rangle} %%% PRODUIT SCALAIRE
\def\f#1,#2,#3{#1\colon#2\to#3} %%%  F:A->B

\def\hfl#1{{\buildrel{#1}\over{\hbox to
12mm{\rightarrowfill}}}}

\title{EXOTIC DEFORMATION QUANTIZATION}

\author{VALENTIN OVSIENKO}

\date{}

\maketitle

\thispagestyle{empty}

\section{Introduction}

Let ${\cal A}$ be one of the following commutative
associative algebras: the algebra of all smooth functions on the plane:
${\cal A}=C^{\infty}({\bf R}^{2})$,
or the algebra of polynomials ${\cal A}={\bf C}[p,q]$ over
${\bf R}$ or ${\bf C}$.
There exists
a non-trivial formal {\it associative} deformation
of ${\cal A}$ called the
{\it Moyal $\star$-product}
(or the standard $\star$-product).
It is defined as an associative
operation ${\cal A}^{\otimes 2}\to {\cal A}[[\hbar]]$
where  $\hbar$ is
a formal variable.
The explicit formula is:
\begin{equation}
F\star_{\hbar}G=
FG+\sum_{k\geq1}\frac{(i\hbar)^k}{2^kk!}\{F,G\}_k,
\label{Moy}
\end{equation}
where
$\{F,G\}_1=
\displaystyle{\frac{\partial F}{\partial p}\frac{\partial G}{\partial q}-
\frac{\partial F}{\partial q}\frac{\partial G}{\partial p}}$
is the standard
Poisson bracket, the higher order
terms are:
\begin{equation}
\{F,G\}_k=
\sum_{i=0}^{k}(-1)^i
\left(
\begin{array}{c}
k\\
i
\end{array}
\right)
\frac{\partial ^kF}{\partial p^{k-i}\partial q^i}
\frac{\partial ^kG}{\partial p^i\partial q^{k-i}}.
\label{New}
\end{equation}

The Moyal product is the {\it unique} (modulo equivalence)
non-trivial
formal deformation of the associative algebra
${\cal A}$ (see \cite{lic}).

\vskip 0,3cm

{\bf Definition 1}. A formal associative deformation
of ${\cal A}$
given by the formula (\ref{Moy}) is called a $\star$-{\it product} if:
\hfill\break
1) the first order term coincides with the Poisson bracket:
$\{F,G\}_1=\{F,G\}$;
\hfill\break
2) the higher order terms $\{F,G\}_k$ are given by differential operators
vanishing on constants: $\{1,G\}_k=\{F,1\}_k=0$;
\hfill\break
3)
$\{F,G\}_k=(-1)^k\{G,F\}_k$.

\vskip 0,3cm

{\bf Definition 2}. Two $\star$-products $\star_{\hbar}$ and $\star_{\hbar}'$
on ${\cal A}$
are called {\it equivalent} if there exists a linear mapping
$A_{\hbar}:{\cal A}\to{\cal A}[[\hbar]]$:
$$
A_{\hbar}(F)=
F+\sum_{k=1}^{\infty}A_k(F)\hbar^k
$$
intertwining the operations $\star_{\hbar}$ and $\star_{\hbar}'$:
$
A_{\hbar}(F)\star_{\hbar}'A_{\hbar}(G)=
A_{\hbar}(F\star_{\hbar}G).
$

\vskip 0,3cm

Consider now
${\cal A}$
as a Lie algebra: the commutator is given by the Poisson bracket.
The Lie algebra ${\cal A}$ has a unique (modulo equivalence)
non-trivial
formal deformation
called the {\it Moyal bracket}
or the Moyal $\star$-commutator:
$
\{F,G\}_t=\frac{1}{i\hbar}(F\star_{\hbar}G-G\star_{\hbar}F),
$
where $t=-\hbar^2/2$.

\vskip 0,3cm

The well-known De Wilde--Lecomte theorem \cite{ld} states
the existence of a non-trivial $\star$-product
for an arbitrary symplectic manifold.
The theory of $\star$-products is a subject of
{\it deformation quantization}.
The geometrical proof of the existence theorem
was given by B. Fedosov \cite{fed}
(see \cite{ew} and \cite{wei} for clear explanation
and survey of recent progress).

\vskip 0,3cm

The
main idea of this paper is to consider the
algebra ${\cal F}(M)$ of functions (with singularities)
on the cotangent bundle
$T^*M$ which are {\it Laurent polynomials} on the fibers.
In contrast to above algebra ${\cal A}$
it turns out that for such algebras
the standard $\star$-product is no more unique
at least if $M$ is one-dimensional: $\dim M=1$.

\vskip 0,3cm

We consider $M=S^1,{\bf R}$ in the real case, and $M=\cal H$
(the upper half-plane)
in the holomorphic case.
The main result of this paper is an explicit construction
of a new $\star$-product on
the algebra ${\cal F}(M)$ non-equivalent
to the standard Moyal product.
This $\star$-product is equivariant with respect to the M\"obius
transformations.
The construction is based on the bilinear
$SL_2$-equivariant operations
on tensor-densities on $M$,
known as {\it Gordan transvectants} and {\it Rankin-Cohen brackets}.

We study the relations between the new $\star$-product and
extensions of the Lie algebra $\Vect(S^1)$.

\vskip 0,5cm

The results of this paper are closely related to those of
the recent work of P. Cohen, Yu. Manin and D. Zagier \cite{cmz}
where a one-parameter family of
associative products on the space of classical modular forms
is constructed using the same $SL_2$-equivariant bilinear operations.

\section{Definition of the exotic $\star$-product}

{\bf 2.1 Algebras of Laurent polynomials}.
Let ${\cal F}$ be one of the following associative algebras
of functions:
$$
{\cal F}=C[p,1/p]\otimes C^{\infty}({\bf R})\;\;\;\;\;
\hbox{or}
\;\;\;\;\;
{\cal F}=C[p,1/p]\otimes \hbox{Hol}({\cal H})
$$
This means, it consists of functions of the type:
\begin{equation}
F(p,q)=\sum_{i=-N}^{N} p^if_i(q)
\label{def}
\end{equation}
where $f_i(q)\in C^{\infty}({\bf R})$
in the real case, or $f_i(q)$
are holomorphic functions on
the upper half-plane ${\cal H}$
or $f_i\in C[q]$ (respectively).

We will also consider the algebra of
polynomials:
${\bf C}[p,1/p,q]$ (Laurent polynomials in $p$) .

\vskip 0,3cm

{\bf 2.2 Transvectants}.
Consider the following bilinear operators
on functions of one variable:
\begin{equation}
J_k^{m,n}(f,g)=
\sum_{i+j=k}(-1)^i
\left(
\begin{array}{c}
k\\
i
\end{array}
\right)
\frac{(2m -i)!(2n -j)!}{(2m -k)!(2n -k)!}
f^{(i)}g^{(j)}
\label{Tra}
\end{equation}
where $f=f(z),g=g(z)$,
$f^{(i)}(z)=\displaystyle\frac{d^if(z)}{dz^i}$.

These operators satisfy a
remarkable property: they are equivariant under
M\"obius (linear-fractional) transformations.
Namely, suppose that the transformation
$z\mapsto\frac{az+b}{cz+d}$ (with $ad-bc=1$)
acts on the arguments as follows:
$$
f(z)\mapsto f\bigg(\frac{az+b}{cz+d}\bigg)(cz+d)^{m},\;\;\;
g(z)\mapsto g\bigg(\frac{az+b}{cz+d}\bigg)(cz+d)^{n},
$$
then $J_k^{m,n}(f,g)$ transforms as:
$J_k^{m,n}(f,g)(z)
\mapsto J_k^{m,n}(f,g)(\frac{az+b}{cz+d})(cz+d)^{m+n-k}.$

In other words,
the operations (\ref{Tra})
are bilinear $SL_2$-equivariant mappings
on tensor-densities:
$$
J^{mn}_k:{\cal F}_m\otimes{\cal F}_n\to{\cal F}_{m+n-k}
$$
where ${\cal F}_l$ is the space of tensor-densities of degree
$-l$: $\phi=\phi(z)(dz)^{-l}$.

\vskip 0,3cm

The operations (\ref{Tra}) were discovered
more than one hundred years ago by Gordan {\cite{Gor} who
called them the {\it transvectants}.
They have been rediscovered many times:
in the theory of modular functions by Rankin \cite{ran}
and by Cohen \cite{coh} (so-called Rankin-Cohen brackets),
in differential projective geometry by Janson and Peetre
\cite{Pee}.
The ``multi-dimensional transvectants''
were defined in \cite{vir} in the context of the
the Virasoro algebra and symplectic and contact geometry.

\vskip 0,3cm

{\bf 2.3 Main definition}.
Define the following bilinear mapping
${\cal F}^{\otimes 2}\to{\cal F}[[\hbar]]$,
for $F=p^mf(q),G=p^ng(q)$, where $m,n\in {\bf Z}$, put:
\begin{equation}
F{\widetilde\star}_{\hbar}G=
\sum_{k=0}^{\infty}
\frac{(i\hbar)^k}{2^{2k}}p^{(m+n-k)}J_k^{m,n}(f,g)
\label{Ovs}
\end{equation}
Note, that the first order term coincides with the Poisson bracket.

This operation will be the main subject of this paper.
We call it the {\it exotic $\star$-product}.

\vskip 0,3cm

{\bf 2.4 Remark}. Another
one-parameter family of operations
on modular forms:
$f\star^{\kappa} g=
\sum_{n=0}^{\infty}t_n^{\kappa}(k,l)J^{kl}_n(f,g),$
where $t_n^{\kappa}(k,l)$ are very interesting and complicated coefficients,
is defined in \cite{cmz}.

\section{Main theorems}
We formulate here the main results of this paper.
All the proofs will be given in Sections 4-7.

\vskip 0,3cm

{\bf 3.1 Non-equivalence}.
The Moyal $\star$-product (\ref{Moy})
defines a non-trivial formal
deformation of ${\cal F}$.
We will show that the formula (\ref{Ovs})
defines a $\star$-product
non-equivalent to the standard Moyal product.

\proclaim Theorem 1. The operation (\ref{Ovs}) is associative, it
defines a formal deformation
of the algebra ${\cal F}$
which is not equivalent to the Moyal product.\par

The associativity of the product (\ref{Ovs})
is a trivial corollary of Proposition 1 below.
To prove the non-equivalence,
we will use the relations with extensions
of the Lie algebra of vector fields on $S^1$:
$\Vect(S^1)\subset{\cal F}$ (cf. Sec.5).

\vskip 0,3cm

It is interesting to note, that the constructed $\star$-product
is equivalent
to the standard Moyal product if we consider it on the algebra
$C^{\infty}(T^*M\setminus M)$ of all smooth functions
(not only Laurent polynomials on fibers), cf. Corollary 1 below.

\vskip 0,3cm

{\bf 3.2 $sl_2$-equivariance}.
The Lie algebra $sl_2({\bf R})$ has two natural embeddings
into the Poisson Lie algebra on ${\bf R}^2$:
the {\it symplectic Lie algebra}
$sp_2({\bf R})\cong sl_2({\bf R})$
generated by quadratic polynomials
$
(p^2,pq,q^2)$
and another one
 with generators:
$(p,pq,pq^2)$ which is called the {\it M\"obius algebra}.

It is well-known that
the Moyal product (\ref{Moy})
is the unique non-trivial formal deformation of the associative algebra
of functions on ${\bf R}^2$ equivariant under the action of the symplectic
 algebra. This means, (\ref{Moy})
satisfies the Leibnitz property:
\begin{equation}
\{F,G\star_{\hbar}H\}=
\{F,G\}\star_{\hbar}H+G\star_{\hbar}\{F,H\}
\label{Lei}
\end{equation}
where $F$ is a quadratic polynomial (note that
$\{F,G\}_t=\{F,G\}$ if $F$ is a quadratic polynomial).

\proclaim Theorem 2. The product (\ref{Ovs}) is the unique formal
deformation of the associative algebra ${\cal F}$
equivariant under the action of the M\"obius algebra.\par

The product (\ref{Ovs}) is the unique non-trivial formal deformation of
${\cal F}$ satisfying (\ref{Lei}) for $F$ from the M\"obius $sl_2$ algebra.

\vskip 0,3cm

{\bf 3.3 Symplectomorphism $\Phi$}.
The relation between the Moyal product and the product (\ref{Ovs}) is
as follows. Consider the symplectic mapping
\begin{equation}
\Phi(p,q)=\bigg(\frac{p^2}{2},\frac{q}{p}\bigg).
\label{Phi}
\end{equation}
defined on
${\bf R}^2\setminus{\bf R}$
in the real case
and on ${\cal H}$ in the complex case.

\proclaim Proposition 1. The product (\ref{Ovs}) is
the $\Phi$-conjugation of the Moyal product:
\begin{equation}
F{\widetilde\star}_{\hbar}G=
F\star_{\hbar}^{\Phi}G:=
(F\circ\Phi\star_{\hbar}G\circ\Phi)\circ\Phi^{-1}
\label{Con}
\end{equation}\par

{\bf Remark}. The mapping (\ref{Phi}) (in the complex case)
can be interpreted as follows. It transforms
the space of {\it holomorphic tensor-densities}
of degree $-k$
on ${\bf CP}^1$
to the space $C^k[p,q]$ of polynomials of degree $k$.
Indeed, there exists a natural isomorphism
$z^n(dz)^{-m}\mapsto p^mq^n$ (where $m\geq2n$)
and $( p^mq^n)\circ\Phi=p^{2m-n}q^n$.

\vskip 0,3cm

{\bf 3.4 Operator formalism}.
The Moyal product is related to the following Weil quantization
procedure.
Define the following differential operators:
\begin{equation}
\matrix{
\widehat p =i\hbar\frac{\partial}{\partial q},\hfill  \cr\noalign{\smallskip}
\widehat q =q\hfill \cr
}
\label{Sch}
\end{equation}
satisfying the canonical relation:
$[\widehat p,\widehat q]=i\hbar \hbox{{\bf I}}$.
Associate to each polynomial $F=F(p,q)$
the differential operator
$
\widehat F=\hbox{Sym}F(\widehat p,\widehat q)
$
symmetric in $\widehat p$ and $\widehat q$.
The Moyal product on the algebra of polynomials
coincides with the product of differential
operators:
$
\widehat {F\star_{\hbar}G}=\widehat F\widehat G.
$

We will show that the $\star$-product (\ref{Ovs}) leads to
the operators:
\begin{equation}
\matrix{
\widehat p^{\Phi} =(\frac{i\hbar}{2})^2\Delta,\hfill  \cr\noalign{\smallskip}
\widehat q^{\Phi} =\frac{1}{4i\hbar}(\Delta^{-1}\circ A+A\circ \Delta^{-1})
\hfill  \cr\noalign{\smallskip}
}
\label{Ope}
\end{equation}
(where $\Delta=\frac{\partial^2}{\partial q^2}$ and
$A=2q\frac{\partial}{\partial q}+1$ is the dilation operator)
also satisfying the canonical relation.

Remark that $\widehat p^{\Phi}$ and $\widehat q^{\Phi}$
given by (\ref{Ope})
on the Hilbert space $L_2({\bf R})$,
are not equivalent
to the operators (\ref{Sch}) since $\widehat q^{\Phi}$ is symmetric but
not self-adjoint (see \cite{dix} on this subject).

\vskip 0,3cm

{\bf 3.5 ``Symplectomorphic'' deformations}.
Let us consider the general situation.

\proclaim Proposition 2. Given a symplectic manifold $V$
endowed with a $\star$-product
$\star_{\hbar}$ and a symplectomorphism $\Psi$ of $V$,
if there exists a
hamiltonian isotopy of $\Psi$
to the identity, then the $\Psi$-conjugate product
$\star_{\hbar}^{\Psi}$ defined according to the formula (\ref{Con})
is equivalent
to $\star_{\hbar}$.\par

\proclaim Corollary 1. The $\star$-product
(\ref{Ovs}) considered on the algebra of all smooth functions
$C^{\infty}(T^*{\bf R}\setminus{\bf R})$
is equivalent to the Moyal product.\par

\vskip 0,3cm

\section{M\"obius-invariance}

In this section we prove Theorem 2.
We show that the operations of transvectant (\ref{Tra})
are $\Phi$-conjugate of the terms of the Moyal product.

\vskip 0,3cm

{\bf 4.1 Lie algebra $\Vect ({\bf R})$ and modules of tensor-densities}.
Let $\Vect ({\bf R})$ be the Lie algebra
of smooth (or polynomial) vector fields on ${\bf R}$:
$$
X=X(x)\frac{d}{dx}
$$
with the commutator
$$
[X(x)\frac{d}{dx},Y(x)\frac{d}{dx}]=
(X(x)Y^{\prime}(x)-X^{\prime}(x)Y(x))\frac{d}{dx}.
$$
The natural embedding of the Lie algebra $sl_2\subset\Vect ({\bf R})$
is generated by the vector fields
$d/dx,xd/dx,x^2d/dx.$

Define a 1-parameter family of $\Vect ({\bf R})$-actions on
$C^{\infty}({\bf R})$ given by
\begin{equation}
L^{(\lambda)}_Xf=X(x)f^{\prime}(x)-\lambda X^{\prime}(x)f(x)
\label{Lie}
\end{equation}
where $\lambda\in {\bf R}$.
Geometrically, $L^{(\lambda)}_X$
is the operator of Lie derivative on {\it tensor-densities}
of degree $-\lambda$:
$$f=f(x)(dx)^{-\lambda}.$$

Denote ${\cal F}_{\lambda}$
the $\Vect ({\bf R})$-module structure on $C^{\infty}({\bf R})$
given by (\ref{Lie}).

\vskip 0,3cm

{\bf 4.2 Transvectant as a bilinear $sl_2$-equivariant operator}.
The operations (\ref{Tra}) can be defined
as bilinear
mappings on $C^{\infty}({\bf R})$ which are
{\it $sl_2$-equivariant}:

\proclaim Statement 4.1. For each $k=0,1,2,\dots$
there exists a unique (up to a constant) bilinear $sl_2$-equivariant
mapping
$$
{\cal F}_{\mu}\otimes{\cal F}_{\nu}\to
{\cal F}_{\mu+\nu-k}.
$$
It is given by $f\otimes g\mapsto J_k^{\mu,\nu}(f,g)$.\par

{\bf Proof}: straightforward (cf. \cite{Gor}, \cite{Pee}).

\vskip 0,3cm

{\bf 4.3 Algebra ${\cal F}$ as a module over $\Vect ({\bf R})$}.
The Lie algebra $\Vect ({\bf R})$ can be considered
as a Lie subalgebra of ${\cal F}$.
The embedding $\Vect ({\bf R})\subset{\cal F}$ is given by:
$$
X(x)\frac{d}{dx}\mapsto pX(q).
$$
The algebra ${\cal F}$
is therefore, a $\Vect ({\bf R})$-module.

\proclaim Lemma 4.2. The algebra ${\cal F}$ is decomposed to
a direct sum of $\Vect ({\bf R})$-modules:
$$
{\cal F}=\oplus_{m\in{\bf Z}}{\cal F}_m.
$$\par

{\bf Proof}. Consider the subspace
of ${\cal F}$ consisting of functions
homogeneous of degree $m$ in $p$:
$F=p^mf(q)$.
This subspace
is a $\Vect ({\bf R})$-module isomorphic to ${\cal F}_m$.
Indeed,
$
\{pX(q),p^mf(q)\}=p^m(Xf^{\prime}-mX^{\prime}f)=p^mL_X^{(m)}f.
$

\vskip 0,3cm

{\bf 4.4 Projective property of the diffeomorphism $\Phi$}.
The transvectants (\ref{Tra}) coincide with
the $\Phi$-conjugate operators
(\ref{New}) from the Moyal product:

\proclaim Proposition 4.3. Let $F=p^mf(q),G=p^ng(q)$, then
\begin{equation}
\Phi^{*-1}\{\Phi^*F,\Phi^*G\}_k=
\frac{k!}{2^k}p^{m+n-k}J_k^{m,n}(f,g).
\label{Iso}
\end{equation}\par

{\bf Proof}.
The symplectomorphism $\Phi$ of ${\bf R}^2$
intertwines the symplectic algebra $sp_2\equiv sl_2$
and the M\"obius algebra:
$
\Phi^*(p,pq,pq^2)=(\frac{1}{2}p^2,\frac{1}{2}pq,\frac{1}{2}q^2).
$
Therefore, the operation $\Phi^{*-1}\{\Phi^*F,\Phi^*G\}_k$
is M\"obius-equivariant.

On the other hand, one has:
$
\Phi^*F=\frac{1}{2^m}p^{2m}f(\frac{q}{p})$
and $\Phi^*G=\frac{1}{2^n}p^{2m}g(\frac{q}{p}).
$
Since $\Phi^*F$ and $\Phi^*G$ are
homogeneous of degree $2m$ and $2n$ (respectively),
the function $\{\Phi^*F,\Phi^*G\}_k$
is also homogeneous of degree $2(m+n-k)$.
Thus, the operation
$\{F,G\}^{\Phi}_k=\Phi^{*-1}\{\Phi^*F,\Phi^*G\}_k$
defines a bilinear mapping on the space of tensor-densities
${\cal F}_m\otimes{\cal F}_n\to{\cal F}_{m+n-k}$
which is $sl_2$-equivariant.

Statement 4.1 implies that it is proportional to $J_k^{m,n}$.
One easily verifies the coefficient of proportionality
for $F=p^m,G=p^nq^k$, to obtain the formula (\ref{Iso}).

Proposition 4.3 is proven.

{\bf Remark}. Proposition 4.3 was proven in \cite{oo}.
We do not know whether this elementary fact has been mentioned
by classics.

\vskip 0,3cm

{\bf 4.5 Proof of Theorem 2}.
Proposition 4.3 implies that the formula (\ref{Ovs}) is a $\Phi$-conjugation
of the Moyal product and is given by the formula (\ref{Con}).

Proposition 1 is proven.

\vskip 0,3cm

It follows that (\ref{Ovs}) is a $\star$-product on ${\cal F}$
equivariant under
the action of the M\"obius $sl_2$ algebra.
Moreover, it is the unique $\star$-product with this property
since the Moyal product is the unique $\star$-product equivariant under
the action of the symplectic algebra.

Theorem 2 is proven.

\section{Relation with extensions of the Lie algebra
$\Vect(S^1)$}

We prove here that the $\star$-product (\ref{Ovs}) is
not equivalent to the Moyal product.

Let $\Vect(S^1)$ be the Lie algebra of vector fields on the circle.
Consider the embedding $\Vect(S^1)\subset{\cal F}$
given by functions on ${\bf R}^2$ of the type:
$X=pX(q)$ where $X(q)$ is periodical: $X(q+1)=X(q)$.

\vskip 0,3cm

{\bf 5.1 An idea of the proof of Theorem 1}.
Consider the formal deformations of the {\it Lie algebra} ${\cal F}$
associated to the $\star$-products (\ref{Moy}) and (\ref{Ovs}).
The restriction of the Moyal bracket
to $\Vect (S^1)$ is identically zero.
We show that the restriction of the $\star$-commutator
$$
\widetilde{\{F,G\}_t}=
\frac{1}{i\hbar}(F{\widetilde\star}_{\hbar}G-G{\widetilde\star}_{\hbar}F),
\;\;\;\;\;
t=-\frac{\hbar^2}{2}
$$
associated to the $\star$-product (\ref{Ovs}) defines a series of non-trivial
extensions of the Lie algebra $\Vect (S^1)$ by the modules
${\cal F}_k(S^1)$ of tensor-densities on $S^1$ of degree $-k$.

\vskip 0,3cm

{\bf 5.2 Extensions and the cohomology group
$H^2(\Vect (S^1);{\cal F}_{\lambda})$}.
Recall that {\it an extension} of a Lie algebra by its module
is defined by a 2-cocycle on it with values in this module.
To define an extension of $\Vect (S^1)$ by the module
${\cal F}_{\lambda}$ one needs therefore, a bilinear mapping
$c:\Vect (S^1)^{\otimes 2}\to{\cal F}_{\lambda}$
which verify the identity $\delta c=0$:
$$
c(X,[Y,Z])+L_X^{(\lambda)}c(Y,Z)+(cycle_{X,Y,Z})=0.
$$
(see \cite{fuc}).

 The cohomology group
$H^2(\Vect (S^1);{\cal F}_{\lambda})$
were calculated in \cite{tsu} (see \cite{fuc}).
This group is trivial for each value of $\lambda$
except $\lambda=0,-1,-2,-5,-7$.
The explicit formul{\ae} for the corresponding non-trivial cocycles are
given in \cite{or2}.
If $\lambda=-5,-7$, then
$\dim H^2(\Vect (S^1);{\cal F}_{\lambda})=1$,
the cohomology group is generated by the unique (up to equivalence)
non-trivial cocycle.
We will obtain these cocycles from the
$\star$-commutator.

\vskip 0,3cm

{\bf 5.3 Non-trivial cocycles on $\Vect (S^1)$}.

Consider the restriction of the $\star$-commutator
$\widetilde{\{\;,\;\}_t}$ (corresponding to
the $\star$-product (\ref{Ovs})) to $\Vect(S^1)\subset{\cal F}$:
let
$$
X=pX(q),\;\;Y=pY(q),
$$
then from (\ref{Ovs}):
$$
\{X,Y\}_t=\{X,Y\}+\sum_{k=1}^{\infty}
\frac{t^k}{2^{2k+1}}\;\frac{1}{p^{2k-1}}\;
J_{2k+1}^{1,1}(X,Y)
$$
It follows from the Jacobi identity that the first non-zero term
of the series $\widetilde{\{X,Y\}_t}$ is a 2-cocycle on $\Vect (S^1)$
with values in one of the $\Vect (S^1)$-modules ${\cal F}_k(S^1)$.

Denote for simplicity $J_{2k+1}^{1,1}$ by $J_{2k+1}$.

{}From the general formula (\ref{Tra}) one obtains:

\proclaim Lemma 5.1. First two terms
of
$\widetilde{\{X,Y\}_t}$
are identically zero: $J_3(X,Y)=0, J_5(X,Y)=0$,
the next two terms are proportional to:
\begin{equation}
\begin{array}{l}
J_7(X,Y) =
X^{\prime\prime\prime}Y^{(IV)}-X^{(IV)}Y^{\prime\prime\prime}\\
J_9(X,Y)=
2(X^{\prime\prime\prime}Y^{(VI)}-X^{(VI)}Y^{\prime\prime\prime})
 -9(X^{(IV)}Y^{(V)}-X^{(V)}Y^{(IV)})
\end{array}
\label{Coc}
\end{equation}\par

The transvectant $J_7$ defines therefore a
2-cocycle. It is a remarkable fact that the same fact is true for $J_9$:

\proclaim Lemma 5.2. (see \cite{or2}). The mappings
$$
J_7:\Vect (S^1)^{\otimes 2}\to{\cal F}_{-5}
\;\;\;\hbox{and}\;\;\;
J_9:\Vect (S^1)^{\otimes 2}\to{\cal F}_{-7}
$$
are 2-cocycles on $\Vect (S^1)$ representing the unique non-trivial
classes of the cohomology groups
$H^2(\Vect (S^1);{\cal F}_{-5})$ and $H^2(\Vect (S^1);{\cal F}_{-7})$
(respectively).\par

{\bf Proof}.
Let us prove that $J_9$ is a 2-cocycle on $\Vect (S^1)$.
The Jacobi identity for the bracket $\{\;,\;\}_t$ implies:
$$
\{X,J_9(Y,Z)\}+J_9(X,\{Y,Z\})+
J_3(X,J_7(Y,Z))\;\;+(cycle_{X,Y,Z})=0
$$
for any $X=pX(q),Y=pY(q),Z=pZ(q)$.
One checks that the expression
$J_3(X,J_7(Y,Z))$ is proportional to
$X'''(Y'''Z^{(IV)}-Y^{(IV)}Z''')$,
and so one gets:
$$
J_3(X,J_7(Y,Z))\;\;+(cycle_{X,Y,Z})=0.
$$
We obtain the
following relation:
$$
\{X,J_9(Y,Z)\}+J_9(X,\{Y,Z\})\;\;+(cycle_{X,Y,Z})=0
$$
which means that $J_9$ is a 2-cocycle. Indeed, recall that for any
tensor density $a$, $\{pX,p^ma\}=p^mL_{Xd/dx}^{(m)}(a)$.
Thus, the last relation coincides
with the relation $\delta J_9=0$.

Let us show now that the cocycle $J_7$ on $\Vect(S^1)$ is not trivial.
Consider a linear differential operator
$A:\Vect(S^1)\to{\cal F}_5$. It is given by:
$A(X(q)d/dq)=(\sum_{i=0}^{K}a_iX^{(i)}(q))(dq)^5$,
then $\delta A(X,Y)=L_X^{(5)}A(Y)-L_Y^{(5)}A(X)-A([X,Y])$.
The higher order part of this expression has a
non-zero term $(5-K)a_KX'Y^{(K)}$ and therefore
$J_7\not=\delta A$.

In the same way one proves that the cocycle $J_9$ on $\Vect(S^1)$ is
non-trivial.

Lemma 5.2 is proven.

It follows that the $\star$-product (\ref{Ovs})
on the algebra ${\cal F}$
is not equivalent to the Moyal product.

Theorem 1 is proven.

\section{Operator representation}

We are looking for an linear mapping (depending on $\hbar$)
$F\mapsto\widehat F^{\Phi}$ of the associative
algebra of Laurent polynomials
${\cal F}=C[p,1/p,q]$ into the algebra of formal
pseudodifferential operators on ${\bf R}$
such that
$$
\widehat {F{\widetilde\star}_{\hbar}G}^{\Phi}=
\widehat F^{\Phi}\widehat G^{\Phi}.
$$

Recall that
the algebra of Laurent polynomials
$C[p,1/p,q]$ with the Moyal product
is isomorphic to the associative algebra
of pseudodifferential operators on ${\bf R}$
with polynomial coefficients (see \cite{adl}).
This isomorphism is
defined on the generators $p\mapsto\widehat p,q\mapsto\widehat q$
by the operators (\ref{Sch}) and $p^{-1}\mapsto\widehat p^{-1}$:
$$
\widehat p^{-1}=
\frac{1}{i\hbar}(\partial/\partial q)^{-1}.
$$

\vskip 0,3cm

{\bf 6.1 Definition}. Put:
\begin{equation}
\widehat F^{\Phi}=\widehat{\Phi^*F}
\label{Qua}
\end{equation}
then $\widehat F^{\Phi}\widehat G^{\Phi}=
\widehat{\Phi^*F}\widehat{\Phi^*G}=\Phi^*F\star_{\hbar}\Phi^*G=
\Phi^*(F{\widetilde\star}_{\hbar}G)=
\widehat {F{\widetilde\star}_{\hbar}G}^{\Phi}$.

One obtains the formul{\ae} (\ref{Ope}). Indeed,
$$
\widehat p^{\Phi}=
\widehat{p^2/2}=
\frac{(i\hbar)^2}{2}\frac{\partial^2}{\partial q^2}.
$$
Since
$
q=
\frac{1}{2}\left((\frac{1}{p}){\widetilde\star}_{\hbar}pq+
pq{\widetilde\star}_{\hbar}(\frac{1}{p})\right),
$
one gets:
$$
\widehat q^{\Phi}=\frac{1}{4i\hbar}(\Delta\circ A+A\circ\Delta)
$$

\vskip 0,3cm

{\bf 6.2 $sl_2$-equivariance}. For the M\"obius $sl_2$ algebra one has:
$$
\begin{array}{lcl}
\widehat p^{\Phi}&=&\displaystyle\frac{(i\hbar)^2}{2}\frac{\partial^2}{\partial
q^2}\\
\widehat{pq}^{\Phi}&=&
\displaystyle\frac{i\hbar}{4}+\frac{i\hbar}{2}q\frac{\partial}{\partial q}\\
\widehat{pq^2}^{\Phi}&=&\displaystyle\frac{q^2}{2}
\end{array}
$$

\proclaim Lemma 6.1. The mapping $F\mapsto\widehat F^{\Phi}$
satisfies the M\"obius-equivariance condition:
$$
\widehat {\{X,F\}}^{\Phi}=
[\widehat X^{\Phi},\widehat F^{\Phi}]
$$
for $X\in sl_2$.\par

{\bf Proof}.
It follows immediately from Theorem 2.
Indeed,
the $\star$-product (\ref{Ovs}) is
$sl_2$-equivariant
(that is, satisfying the relation:
$\{X,F\}_t=\{X,F\}$ for $X\in sl_2$).

\vskip 0,3cm

{\bf Remark}. Beautiful explicit formul{\ae} for
$sl_2$-equivariant mappings from the space of
tensor-densities to the space of pseudodifferential
operators are given in \cite{cmz}.

\section{Hamiltonian isotopy}
The simple calculations below are quite standard for
the cohomological technique. We need them to prove
Corollary 1 (from Sec. 3).

Given a syplectomorphism $\Psi$ of a symplectic manifold $V$
and a formal
deformation $\{\;,\;\}_t$ of the Poisson bracket on $V$,
we prove that if $\Psi$ is isotopic to the identity then
the formal
deformation $\{\;,\;\}_t^{\Psi}$ defined by:
$$
\{F,G\}_t^{\Psi}=\Psi^{*-1}\{\Psi^*F,\Psi^*G\}_t
$$
is equivalent to $\{\;,\;\}_t$. The similar proof is valid
in the case of
$\star$-products.

 Recall that two symplectomorphisms $\Psi$ and $\Psi'$
of a symplectic manifold $V$
are {\it isotopic}
if there exists a family of functions $H_{(s)}$ on $V$ such that
the symplectomorphism $\Psi_1\circ\Psi_2^{-1}$
is the flow of the Hamiltionian vector field with the Hamiltonian function
$H_{(s)}, 0\leq s\leq 1$.

 Let $\Psi_{(s)}$ be
the the flow of a family of functions $H=H_{(s)}$.
We will prove that the equivalence class of the formal
deformation $\{\;,\;\}_t^{\Psi_{(s)}}$ does not depend on $s$.

\vskip 0,3cm

{\bf 7.1 Equivalence of homotopic cocycles}.
Let us first show that the cohomology class
of the cocycle
$C_3^{\Psi_{(s)}}$:
$$
C_3^{\Psi_{(s)}}(F,G)=\Psi^{*-1}_{(s)}C_3(\Psi^*_{(s)}F,\Psi^*_{(s)}G)
$$
 does not depend on $s$.
To do this, it is sufficient to prove
that the derivative $\dot C_3=\frac{d}{ds}C_3^{\Psi_{(s)}}|_{s=0}$
is a coboundary. One has
$$
\dot C_3(F,G)=
C_3(\{H,F\},G)+C_3(F,\{H,G\})-\{H,C_3(F,G)\}.
$$
The relation $\delta C_3=0$ implies:
$$
\dot C_3(F,G)=\{F,C_3(G,H)\}-\{G,C_3(F,H)\}-C_3(\{F,G\},H)
$$
This means,
$
\frac{d}{ds}C_3^{\Psi_{(s)}}|_{s=0}=
\delta B_{H}
$,
where $B_{H}(F)=C_3(F,H)$.

\vskip 0,3cm

{\bf 7.2 General case}.
Let us apply the same arguments to prove that the deformations
$\{\;,\;\}_t^{\Psi_{(s)}}$ are equivalent to each other for all values
of $s$. We must show that there exists a family of mappings
$
A_{(s)}(F)=F+\sum_{k=1}^{\infty}{A_{(s)}}_k(F)t^k
$
such that $A_{(s)}^{-1}(\{A_{(s)}(F),A_{(s)}(G)\}_t)=\{F,G\}_t$.

It is sufficient to verify the existence of a mapping
$
a(F)=\sum_{k=1}^{\infty}a_k(F)t^k
$
(the derivative: $a(F)=d/ds(A_{(s)}(F))|_{s=s_0}$)
such that
$$
\frac{d}{ds}\{F,G\}_t^{\Psi_{(s)}}|_{s=s_0}=\{a(F),G\}_t+\{F,a(G)\}-a(\{F,G\}_t
$$
One has:
$$
\frac{d}{ds}\{F,G\}_t^{\Psi_{(s)}}|_{s=s_0}=\{\{F,H\},G\}_t+\{F,\{G,H\}\}_t-\{\{F,G\}_t,H\}
$$
{}From the Jacobi identity:
$$
\{\{F,H\}_t,G\}_t+\{F,\{G,H\}_t\}_t-\{\{F,G\}_t,H\}_t=0
$$
one obtains that the mapping $a(F)$ can be written in the form:
$$
a(F)=\sum_{k=1}^{\infty}\frac{1}{(2k+1)!}C_{2k+1}(F,H_{(s_0)})t^k
$$

\vskip 0,3cm

{\bf 7.3 Proof of Corollary 2}.
Consider the $\star$-product (\ref{Con}) given by
$F\star_{\hbar}^{\Phi}G$,
where $F\star_{\hbar}G$ is the Moyal product (\ref{Moy}),
$\Phi:(p,q)\mapsto(p^2/2,q/p)$.
It is defined on ${\bf R}^2\setminus{\bf R}$.

The $\star$-product
(\ref{Con}) on the algebra
$C^{\infty}({\bf R}^2\setminus{\bf R})$
is equivalent to the Moyal product.
Indeed, the symplectomorphism $\Phi$
is isotopic to the identity
in the group of all smooth symplectomorphisms of ${\bf R}^2\setminus{\bf R}$.
The isotopy is: $\Phi_s:(p,q)\mapsto(\frac{p^{1+s}}{1+s},\frac{q}{p^s})$,
where $s\in[0,1]$.

Recall, that the $\star$-product $F\star_{\hbar}^{\Phi}G$
on the algebra ${\cal F}$ is not equivalent to
the Moyal product since it coincides with the product (\ref{Ovs}).

The family $\Phi_s$ does not preserve
the algebra ${\cal F}$. Theorem 1 implies that $F$ is not
isotopic to the identity
in the group of symplectomorphisms of ${\bf R}^2\setminus{\bf R}$
preserving the algebra ${\cal F}$.

\section{Discussion}

{\bf 8.1 Difficulties in multi-dimensional case}

There exist multi-dimensional analogues of transvectants
\cite{vir} and \cite{or1}.

Consider the projective space ${\bf RP}^{2n+1}$
endowed with the standard contact structure
(or an open domain of the complex projective space ${\bf CP}^{2n+1}$).
There exists an unique bilinear differential operator of order
$k$ on tensor-densities
equivariant with respect to the action of the group $Sp_{2n}$
(see \cite{vir}, \cite{or1}):
\begin{equation}
J_k^{\lambda,\mu}:
{\cal F}_{\lambda}\otimes{\cal F}_{\mu}\to
{\cal F}_{\lambda+\mu-\frac{k}{n+1}}
\label{mul}
\end{equation}
where
${\cal F}_{\lambda}={\cal F}_{\lambda}({\bf P}^{2n+1})$
is the space of tensor-densities on
${\bf P}^{2n+1}$ of degree $-\lambda$:
$$
f=f(x_1,\dots,x_{2n+1})(dx_1\wedge\dots dx_{2n+1})^{-\lambda}.
$$

The space of tensor-densities ${\cal F}_{\lambda}({\bf RP}^{2n+1})$
is {\it isomorphic as a module over
the group of contact diffeomorphisms}
to the space of homogeneous functions
on ${\bf R}^{2n+2}$. The isomorphism is given by:
$$
f\mapsto F(y_1,\dots,y_{2n+2})=
y_{2n+2}^{-\lambda(n+1)}f(\frac{y_1}{y_{2n+2}},
\dots,\frac{y_{2n+1}}{y_{2n+2}})
$$
Then, the operations (\ref{mul}) are defined as the
restrictions of the terms of the standard $\star$-product on ${\bf R}^{2n+2}$.

\vskip 0,3cm

The same formula (\ref{Ovs}) defines a
$\star$-product on the space of tensor-densities on ${\bf CP}^{2n+1}$)
(cf. \cite{or1}).
However, there is no analogues of the symplectomorphism (\ref{Phi}).
I do not know if there exists
a $\star$-product on the Poisson algebra
$C[y_{2n+2},y_{2n+2}^{-1}]\otimes C^{\infty}({\bf RP}^{2n+1})$
non-equivalent to the standard.

\vskip 0,3cm

{\bf 8.2 Classification problem}.
The classification (modulo equivalence)
of $\star$-products on the
Poisson algebra ${\cal F}$ is an interesting open problem.
It is related to calculation
of cohomology groups
$H^2({\cal F};{\cal F})$ and $H^3({\cal F};{\cal F})$.
The following result is announced in \cite{dzh}:
$\dim H^2({\cal F};{\cal F})=2$.

Let us formulate a conjecture in the compact case.
Consider the Poisson algebra
${\cal F}(S^1)$ of functions on $T^*S^1\setminus S^1$
which are Laurent polynomials on the fiber:
$F(p,q)=\sum_{-N\leq i\leq N}p^if_i(q)$ where
$f_i(q+1)=f_i(q)$.
\proclaim Conjecture. Every $\star$-product
on ${\cal F}(S^1)$ is equivalent to (\ref{Moy})
or to (\ref{Ovs}).\par

\vskip 0,5cm

I am grateful to Yu. I. Manin
who explained me the notion of Rankin-Cohen
brackets in the theory of modular forms,
for important remarks, references and
clarifying discussions.
It is a pleasure to thank
O. Ogievetsky and C. Roger
for multiple stimulation collaboration;
M. Audin,
C. Duval, P. Lecomte, E. Mourre, S. Tabachnikov,
P. Seibt, Ya. Soibel'man and F. Ziegler
for fruitful discussions
and
their interest to this work.

\vskip 1cm

%%%%%%%%%%%%%%%%%%%%%%%%%%%%%%%%%%%%%%%%%%%%%%%%%%%%%%%%%%%%%%%%%%%%%%%%%%%%%%
%%%%%%%%%%%%%%%%%%%%%%%%%%%%%%%%%%%%%%%%%%%%%%%%%%%%%%%%%%%%%%%%%%%%%%%%%%%%%%

{\sc Centre National De La Recherche Scientifique, Marseille}

\end{document}